\documentclass[10pt,a4paper]{article}
\usepackage{amsmath,amssymb,amsthm}
\usepackage{indentfirst}

\title{Quasi-Invariant Optimal Control Problems\footnote{Partially
  presented at the Invited Session ``Control, Optimization and Computation'',
  10th Mediterranean Conference on Control and Automation (MED2002), Lisbon,
  July 9--12, 2002. Accepted for publication (21.Feb.03) in the journal
  \emph{Portugali\ae  Mathematica}.}}
\author{Delfim F.~M.~Torres \\
         \href{mailto:delfim@mat.ua.pt}{\texttt{delfim@mat.ua.pt}}}
\date{Departamento de Matem\'{a}tica \\
      Universidade de Aveiro \\
      3810-193 Aveiro, Portugal}

\usepackage[dvips]{hyperref}
\hypersetup{colorlinks=false,hyperindex,plainpages=false,
  pdftitle={Quasi-Invariant Optimal Control Problems},
  pdfsubject={Partially presented at the Invited Session
         ``Control, Optimization and Computation'',
         10th Mediterranean Conference on Control
         and Automation (MED2002), Lisbon,
         July 9--12, 2002.},
  pdfkeywords={optimal control, Pontryagin extremals, Noether theorem,
               conservation laws, invariance up to first-order terms in the parameters.},
  pdfauthor={Delfim F.~M.~Torres},
  pdfcreator={\LaTeX\ with package \flqq hyperref\frqq}
}

\newtheorem{theorem}{Theorem}[section]

\newtheorem{proposition}[theorem]{Proposition}

\theoremstyle{remark}
  \newtheorem{remark}{Remark}[section]

\theoremstyle{definition}
  \newtheorem{definition}{Definition}[section]
  \newtheorem{example}{Example}[section]

\newcommand{\ds}{\displaystyle}

\newcommand\MR[1]{\href{http://www.ams.org/mathscinet-getitem?mr=#1}{\textsf{MR}~#1}}
\newcommand\ZBL[1]{\href{http://www.emis.de/MATH-item?#1}{\textsf{Zbl}~#1}}
\newcommand\JFM[1]{\href{http://www.emis.de/cgi-bin/JFM-item?#1}{\textsf{JFM}~#1}}

\begin{document}

\maketitle


\begin{abstract}
We study in optimal control the important relation between invariance
of the problem under a family of transformations, and the existence
of preserved quantities along the Pontryagin extremals.
Several extensions of Noether theorem are provided,
in the direction which enlarges the scope of its application.
We formulate a more general version of Noether's theorem
for optimal control problems, which
incorporates the possibility to consider a
family of transformations depending on several parameters and,
what is more important, to deal with quasi-invariant and
not necessarily invariant optimal control problems. We trust
that this latter extension provides new possibilities and
we illustrate it with several examples, not covered by the previous
known optimal control versions of Noether's theorem.
\end{abstract}







\section{Introduction}

The study of invariant variational problems
\begin{equation*}
J\left[x(\cdot)\right] = \int_a^b L\left(t,x(t),\dot{x}(t)\right)
\mathrm{d}t \longrightarrow \min
\end{equation*}
in the calculus of variations
was initiated in 1918 by Emmy Noether who,
influenced by the works of Klein and Lie on the transformation properties
of differential equations, published
in her gorgeous paper \cite{JFM46.0770.01,MR53:10538} a
fundamental, and now classical result, known as \emph{Noether's theorem}.
The universal principle described by Noether's theorem
(see \textrm{e.g.} \cite[pp.~262--266]{MR16:426a},
\cite[\S 4.3.]{MR87m:49004}, or \cite[\S 20]{ZBL0964.49001}), asserts that
invariance of the integral functionals of the calculus of variations
with respect to a family of transformations result in
existence of a certain conservation
law or equivalently a first integral of the corresponding
Euler-Lagrange differential equations.
This means that the invariance hypothesis leads to quantities,
computed in terms of the Lagrangian and the family of transformations,
which are constant along the extremals.
This result is of great importance in physics, engineering,
systems and control and their applications
(see \cite{MR53:10537,MR58:18024,MR2001a:35105,MR2001k:93034}).
One important application of the Noether theorem is, for example, to
the $n$-body problem. For a discussion of this problem, and
interpretation of the respective first integrals from invariance
under Galilean transformations and application of Noether's
theorem, we refer the reader to \cite{MR1822958} and
\cite[pp.~190--192]{MR98b:49002a} or \cite[Ch.~2]{MR58:18024}.

In the optimal control setting, the relation between invariance of a problem
and the existence of expressions which are constant along any of its extremals,
has been obtained in the publications
by van der Schaft \cite{MR83k:49054} and Sussmann \cite{MR96i:49037},
following the classical Noether's approach based on
the transversality conditions\footnote{In the calculus of variations,
transversality conditions are expressed by the so called
\emph{general variation of the functional}
(see e.g. \cite[\S 13]{ZBL0964.49001} or \cite[p.~185]{MR98b:49002a}).}
(\textrm{cf.} \cite{MR1806135,MR1823009}).
Using the original paper of Emmy Noether \cite{JFM46.0770.01,MR53:10538}
and the more simpler and direct approach of Andrzej Trautman \cite{MR36:3530},
Hanno Rund \cite{MR46:9838} (see also \cite{MR57:13703}) and
John David Logan \cite{MR58:18024} for insight and motivation,
extensions to the previous known optimal control
versions of Noether's theorem were obtained by the present author
in \cite{delfimNoether,delfim3NCNW,delfimEJC}.
Here we attempt to enlarge the range of application of the theorems,
extending the very concept of invariance (Definition~\ref{d:invariant})
by allowing several parameters and equalities up to first-order terms in the
parameters (quasi-invariance).
This extension allows one to formulate a Noether type theorem for
optimal control problems (Theorem~\ref{r:consLaw}) in a much broader
way, enlarging the scope of its application.
Examples not covered by the previous optimal control versions
of Noether's theorem are provided in detail.


\section{The Maximum Principle}

Consider the following optimal control problem,
denoted in the sequel by $(P)$ :
to minimize the integral functional
\begin{equation*}
J\left[x(\cdot),u(\cdot)\right] =
\int_a^b L\left(t,x(t),u(t)\right) dt
\end{equation*}
over the class $W_{1,1}^{n}$ of absolutely continuous
state trajectories
$x(\cdot) = \left(x_1(\cdot),\ldots,x_n(\cdot)\right)$
mapping $[a,b]$ to $\mathbb{R}^n$, and the class
$L_{\infty}^m$ of measurable and essentially bounded controls
$u(\cdot) = \left(u_1(\cdot),\ldots,u_m(\cdot)\right)$
mapping $[a,b]$ to a given set $\Omega \subseteq \mathbb{R}^m$,
subject to the dynamic control system
\begin{equation*}
\dot{x}(t) = \varphi\left(t,x(t),u(t)\right)\quad \text{ for a.a. }
t \in [a,b]\, ,
\end{equation*}
where $L$ and $\varphi$ are assumed to be $C^1$.

The next theorem gives a summary of the celebrated Pontryagin maximum principle
\cite{MR29:3316b}, which is the first-order necessary optimality
condition of optimal control theory.

\begin{theorem}[Pontryagin maximum principle]
\label{th:PMP}
Let $\left(x(\cdot),u(\cdot)\right)$ be a minimizer of the optimal
control problem $(P)$. Then, there exists a nonzero pair
$\left(\psi_0,\psi(\cdot)\right)$, where $\psi_0 \le 0$ is a constant
and $\psi(\cdot)$ a $n$-vector absolutely continuous function with
domain $[a,b]$, such that the following hold for almost all $t$ on
the interval $[a,b]$:
\begin{description}
\item[(i)] the Hamiltonian system
\begin{eqnarray*}
\left\{
\begin{array}{lcl}
\dot{x}(t) & = & \displaystyle
\frac{\partial H\left(t,x(t),u(t),\psi_0,\psi(t)\right)}{\partial
\psi} \, , \\[0.15in]
\dot{\psi}(t) & = & - \displaystyle
\frac{\partial H\left(t,x(t),u(t),\psi_0,\psi(t)\right)}{\partial x}
\, ;
\end{array}
\right.
\end{eqnarray*}
\item[(ii)] the maximality condition
\begin{eqnarray*}
H\left(t,x(t),u(t),\psi_0,\psi(t)\right) = \max_{v \in \Omega}
H\left(t,x(t),v,\psi_0,\psi(t)\right) \, ;
\end{eqnarray*}
\end{description}
with the Hamiltonian
$H(t,x,u,\psi_0,\psi) = \psi_0 L(t,x,u) + \psi \cdot \varphi(t,x,u)$.
\end{theorem}

\begin{definition}
A quadruple $\left(x(\cdot),\,u(\cdot),\,\psi_{0},\,\psi(\cdot)\right)$
satisfying the Hamiltonian system and
the maximality condition is called a (Pontryagin) extremal.
\end{definition}

\begin{remark}
Depending on the specific boundary conditions under consideration in problem $(P)$,
transversality conditions may also appear in the Pontryagin maximum principle.
As far as the results obtained are valid for arbitrary boundary conditions
and the methods which we will employ do not require the use of such transversality
conditions, they are not included in Theorem~\ref{th:PMP}.
\end{remark}


\section{The Quasi-Invariance Definition}
\label{s:preservedQ}

The following notion generalizes the invariance definitions used
in previous versions of Noether's theorem up to first-order terms
in the $r$ parameters $s_1, \cdots, s_r$
(\textrm{cf. e.g.} \cite[Definition 5]{delfim3NCNW}).

\begin{definition}
\label{d:invariant}
If there exists a $C^1$ smooth $r$-parameter family of transformations
\begin{gather}
h^s : [a,b] \times \mathbb{R}^n \times \Omega \rightarrow
       \mathbb{R} \times \mathbb{R}^n \times \mathbb{R}^m \, , \notag \\
h^s(t,x,u) = \left(T(t,x,u,s), X(t,x,u,s), U(t,x,u,s)\right) \, , \label{familyT} \\
s = \left(s_1, \ldots, s_r\right) \, ,
\left\|s\right\| = \sqrt{\sum_{k=1}^{r} \left(s_k\right)^2} < \varepsilon \, , \notag
\end{gather}
which for $s = 0$ reduce to the identity map, $h^0(t,x,u) = (t,x,u)$ for all
$(t,x,u) \in [a,b] \times \mathbb{R}^n \times \Omega$, and
satisfying
\begin{gather}
\begin{split}
L\left(t,x(t),u(t)\right) + \frac{d}{dt} F&\left(t,x(t),u(t),s\right) + o(s) \\
&= L \circ h^s\left(t,x(t),u(t)\right) \frac{d}{dt} T\left(t,x(t),u(t),s\right) \, ,
\end{split} \label{eq:invi} \\
\frac{d}{dt} X\left(t,x(t),u(t),s\right) + o(s)
= \varphi \circ h^s\left(t,x(t),u(t)\right)
\frac{d}{dt} T\left(t,x(t),u(t),s\right) \, , \label{eq:invii}
\end{gather}
for some function $F$ of class $C^1$ and where $o(s)$ denote terms
which go to zero faster than $\left\|s\right\|$, i.e.,
\begin{equation}
\label{eq:limos0}
\lim_{\left\|s\right\|\rightarrow 0} \frac{o(s)}{\left\|s\right\|} = 0 \, ,
\end{equation}
then problem $(P)$ is said to be quasi-invariant under $h^s$.
\end{definition}

\begin{remark}
The types of invariance transformations that we consider are
transformations of the $\left(t,x_1,\ldots,x_n,u_1,\ldots,u_m\right)$-space
which depend upon $r$ small real independent parameters $s_1,\ldots,s_r$.
In Noether's original paper \cite{JFM46.0770.01,MR53:10538}, as well as in
more recent treatments of invariant problems of optimal control
(\textrm{e.g.} \cite{MR1806135}), it is assumed that the transformations
form a \textit{group}.
In the present work, however, we follow the approaches in \cite{delfimEJC}
and \cite{delfim3NCNW}
and we make less stringent assumptions on the transformations -- the group
concept is not required for the investigation of quasi-invariant optimal
control problems.
\end{remark}

The following example shows an optimal control problem quasi-invariant under a
one-parameter family of transformations, in the sense of
Definition~\ref{d:invariant}, but not invariant under all previous
invariance definitions
\cite{MR83k:49054,MR96i:49037,MR1806135,MR1823009,delfimNoether,delfim3NCNW,delfimEJC}
used in connection with the Noether theorem.
This is due to the fact that the integral is not invariant,
but rather invariant up to an exact differential and to first-order
terms in the parameter $s$; while the third component $\varphi_3$
of the phase velocity vector is also invariant only up to first-order
terms in the parameter (quasi-invariant).
\begin{example}[$n = 3$, $m = 2$]
\label{ex:invFO}
We consider problem $(P)$ with $L=u_1^2 + u_2^2$ and
$\varphi = \left(u_1,u_2,\frac{u_2 x_2^2}{2}\right)^T$:
\begin{gather*}
\int_a^b \left(u_1(t)\right)^2 + \left(u_2(t)\right)^2 dt
\longrightarrow \min \, , \\
\begin{cases}
\dot{x_1}(t) = u_1(t) \, , \\
\dot{x_2}(t) = u_2(t)  \, , \\
\dot{x_3}(t) = \displaystyle \frac{u_2(t) \left(x_2(t)\right)^2}{2} \, .
\end{cases}
\end{gather*}
Direct calculations show that the problem is invariant under
$h^s(t,x_1,x_2,x_3,u_1,u_2)
= \left(t,x_1+st,x_2+st,x_3+\frac{1}{2} x_2^2 s t,u_1+s,u_2+s\right)$:
\begin{gather*}
h^0(t,x_1,x_2,x_3,u_1,u_2) = (t,x_1,x_2,x_3,u_1,u_2) \, , \\
L \circ h^s \, \frac{d}{dt}(t) = \left(u_1+s\right)^2+\left(u_2+s\right)^2
= \left(u_1^2 + u_2^2\right)  +
2s\left(u_1 + u_2\right) + 2 s^2 \, ,
\end{gather*}
and equation \eqref{eq:invi} is satisfied with
$F(x_1,x_2,s) = 2s\left(x_1 + x_2\right)$ and
$o(s) = 2 s^2$;
\begin{gather*}
\begin{split}
\varphi_1 \circ h^s \, \frac{d}{dt}(t) &= u_1 + s = \frac{d}{dt} (x_1 + st) \, , \\
\varphi_2 \circ h^s \, \frac{d}{dt}(t) &= u_2 + s = \frac{d}{dt} (x_2 + st) \, , \\
\varphi_3 \circ h^s \, \frac{d}{dt}(t) & =
\frac{\left(u_2 + s\right)\left(x_2+st\right)^2}{2} \\
&= \frac{u_2 x_2^2}{2} + \frac{1}{2} s \left(x_2^2 + 2 x_2 u_2 t\right)
+ \frac{\left(u_2 t^2 + 2 x_2 t\right) s^2 + t^2 s^3}{2} \\
&= \frac{d}{dt}\left(x_3+\frac{1}{2} x_2^2 s t\right) + o(s) \, ,
\end{split}
\end{gather*}
$\left(o(s) = \frac{\left(u_2 t^2 + 2 x_2 t\right) s^2 + t^2 s^3}{2}\right)$
and \eqref{eq:invii} is also satisfied.
\end{example}


\section{The Fundamental Invariance Theorem}

The next fundamental theorem is useful in many ways:
to derive conservation laws for a given quasi-invariant problem $(P)$
(we will see in Section~\ref{sec:NT} how Theorem~\ref{r:fundProp} provide a simple
and direct access  to a Noether theorem -- Theorem~\ref{r:consLaw}) and to
give conditions which allow us  to determine a
family of transformations under which a given optimal
control problem is quasi-invariant
(see Examples~\ref{exn4r2} and \ref{exn4r2PrbNo}, and the ones in
Section~\ref{s:IllExmp}).
If only the transformations are known, equations
\eqref{eq:ds0i} and \eqref{eq:ds0ii} represent first-order
partial differential equations in the unknown functions
$L$ and $\varphi$, and the fundamental theorem can be used to characterize
a set of optimal control problems which possess given invariance properties
(\textrm{cf.} \cite[\S 4.2]{torresCM01I14E}).

\begin{theorem}
\label{r:fundProp}
Necessary conditions for problem $(P)$ to be quasi-invariant under the $r$-parameter
family of transformations \eqref{familyT} are ($k=1,\ldots,r$)~:
\begin{gather}
\frac{d}{dt} \left.\frac{\partial F}{\partial s_k}\right|_{s = 0}
=\frac{\partial L}{\partial t}
      \left.\frac{\partial T}{\partial s_k}\right|_{s = 0}
      + \frac{\partial L}{\partial x} \cdot
      \left.\frac{\partial X}{\partial s_k}\right|_{s = 0}
      + \frac{\partial L}{\partial u} \cdot
      \left.\frac{\partial U}{\partial s_k}\right|_{s = 0}
      + L \frac{d}{dt} \left.\frac{\partial T}{\partial s_k}\right|_{s =
      0}\, , \label{eq:ds0i} \\
\frac{d}{dt} \left.\frac{\partial X}{\partial s_k}\right|_{s = 0}
= \frac{\partial \varphi}{\partial t}
      \left.\frac{\partial T}{\partial s_k}\right|_{s = 0}
      + \frac{\partial \varphi}{\partial x} \cdot
      \left.\frac{\partial X}{\partial s_k}\right|_{s = 0}
      + \frac{\partial \varphi}{\partial u} \cdot
      \left.\frac{\partial U}{\partial s_k}\right|_{s = 0}
      + \varphi \frac{d}{dt} \left.\frac{\partial T}{\partial s_k}\right|_{s =
      0} \, . \label{eq:ds0ii}
\end{gather}
\end{theorem}

\begin{proof}
The proof follows as a simple exercise from the definition of quasi-invariance:
using $h^0(t,x,u) = (t,x,u)$, it suffices to differentiate \eqref{eq:invi}
and \eqref{eq:invii} with respect to $s_k$ and then set $s = 0$.
\end{proof}

\begin{remark}
We are assuming in Theorems~\ref{r:fundProp} and \ref{r:consLaw}
the possibility to reverse the order of differentiation.
\end{remark}

\begin{remark}
From \eqref{eq:invi} one has
\begin{equation*}
o(s) =
L \circ h^s \, \frac{d}{dt} T\left(t,x(t),u(t),s\right)
- L - \frac{d}{dt} F\left(t,x(t),u(t),s\right)\, ,
\end{equation*}
while from \eqref{eq:invii} one obtains
\begin{equation*}
o(s) = \varphi \circ h^s \,
\frac{d}{dt} T\left(t,x(t),u(t),s\right)
- \frac{d}{dt} X\left(t,x(t),u(t),s\right) \, .
\end{equation*}
From these equalities, explicit formulas for the derivatives
of each $o(s_1,\ldots,s_r)$ with respect to $s_k$ ($k=1,\ldots,r$)
can be found. The derivatives vanish for
$s = (s_1,\ldots,s_r) = 0$ due to \eqref{eq:limos0}.
\end{remark}

The next two examples illustrate how Theorem~\ref{r:fundProp}
can be used to guess a family of transformations which maintain
the problem invariant in the sense of Definition~\ref{d:invariant}.
Once again, the possibility of invariance up to first-order terms
in the parameter (quasi-invariance) is crucial.

\begin{example}
\label{exn4r2}
($n = 4$, $m = 2$)
Let us consider the problem
\begin{gather*}
\int_a^b \left(\left(u_1(t)\right)^2 + \left(u_2(t)\right)^2 \right) \, dt
\longrightarrow \min \, , \\
\begin{cases}
\dot{x_1}(t) = x_3(t) \\
\dot{x_2}(t) = x_4(t) \\
\dot{x_3}(t) = - x_1(t) \left(\left(x_1(t)\right)^2 +
\left(x_2(t)\right)^2\right) + u_1(t) \\
\dot{x_4}(t) = - x_2(t) \left(\left(x_1(t)\right)^2 +
\left(x_2(t)\right)^2\right) + u_2(t) \, ,
\end{cases}
\end{gather*}
and look for a one-parameter family of transformations without
changing the time-variable ($T = t$) and with $F \equiv 0$,
under which the problem is quasi-invariant. Theorem~\ref{r:fundProp}
asserts that the following conditions must hold:
\begin{equation}
\label{e:condNecExn4r2}
\begin{cases}
\ds u_1 \left.\frac{\partial U_1}{\partial s}\right|_{s = 0}
     &= \ds - u_2 \left.\frac{\partial U_2}{\partial s}\right|_{s = 0} \\[0.3cm]
\ds \frac{d}{dt} \left.\frac{\partial X_1}{\partial s}\right|_{s = 0}
     &= \ds  \left.\frac{\partial X_3}{\partial s}\right|_{s = 0} \\[0.3cm]
\ds \frac{d}{dt} \left.\frac{\partial X_2}{\partial s}\right|_{s = 0}
     &= \ds  \left.\frac{\partial X_4}{\partial s}\right|_{s = 0} \\[0.3cm]
\ds \frac{d}{dt} \left.\frac{\partial X_3}{\partial s}\right|_{s = 0}
     &= \ds -\left(3 x_1^2 + x_2^2\right)
     \left.\frac{\partial X_1}{\partial s}\right|_{s = 0}
     - 2 x_1 x_2 \left.\frac{\partial X_2}{\partial s}\right|_{s = 0}
     + \left.\frac{\partial U_1}{\partial s}\right|_{s = 0} \\[0.3cm]
\ds \frac{d}{dt} \left.\frac{\partial X_4}{\partial s}\right|_{s = 0}
     &= \ds - 2 x_1 x_2 \left.\frac{\partial X_1}{\partial s}\right|_{s = 0}
      -\left(x_1^2 + 3 x_2^2\right)
      \left.\frac{\partial X_2}{\partial s}\right|_{s = 0}
     + \left.\frac{\partial U_2}{\partial s}\right|_{s = 0} \, .
\end{cases}
\end{equation}
One easily obtains that \eqref{e:condNecExn4r2} is satisfied for
\begin{gather*}
\left.\frac{\partial U_1}{\partial s}\right|_{s = 0} = - u_2 \, , \quad
\left.\frac{\partial U_2}{\partial s}\right|_{s = 0} = u_1 \, , \\
\left.\frac{\partial X_1}{\partial s}\right|_{s = 0} = - x_2 \, , \quad
\left.\frac{\partial X_2}{\partial s}\right|_{s = 0} = x_1 \, , \quad
\left.\frac{\partial X_3}{\partial s}\right|_{s = 0} = - x_4 \, , \quad
\left.\frac{\partial X_4}{\partial s}\right|_{s = 0} = x_3 \, .
\end{gather*}
Choosing $U_1 = u_1 - u_2 s$, $U_2 = u_2 + u_1 s$,
$X_1 = x_1 - x_2 s$, $X_2 = x_2 + x_1 s$, $X_3 = x_3 - x_4 s$,
$X_4 = x_4 + x_3 s$, one can verify that conditions \eqref{eq:invi}
and \eqref{eq:invii} are indeed true:
\begin{gather*}
L \circ h^s \, \frac{d}{dt} T  = \left(u_1 - u_2 s\right)^2 +
\left(u_2 + u_1 s\right)^2 = \left(u_1^2 + u_2^2\right)  +
\left(u_1^2 + u_2^2\right) s^2 = L + o(s) \, , \\
\varphi_1 \circ h^s \, \frac{d}{dt} T =
x_3 - x_4 s = \frac{d}{dt} \left(x_1 - x_2 s\right) = \frac{d}{dt} X_1 \, , \\
\varphi_2 \circ h^s \, \frac{d}{dt} T = x_4 + x_3 s
= \frac{d}{dt} (x_2 + x_1 s) = \frac{d}{dt} X_2 \, , \\
\begin{split}
\varphi_3 \circ h^s \, \frac{d}{dt} T &=
- \left(x_1 - x_2 s\right)\left( (x_1 - x_2 s)^2 + (x_2 + x_1 s)^2\right) + u_1 - u_2 s\\
&= \left[- x_1 (x_1^2 + x_2^2) + u_1 + x_2 (x_1^2 + x_2^2) s - u_2 s\right]
+ \left[ (x_2 s - x_1) (x_1^2 + x_2^2) s^2\right] \\
&= \frac{d}{dt} X_3 + o(s) \, ,
\end{split} \\
\begin{split}
\varphi_4 \circ h^s \, \frac{d}{dt} T &=
- \left(x_2 + x_1 s\right)\left( (x_1 - x_2 s)^2 + (x_2 + x_1 s)^2\right) + u_2 + u_1 s\\
&= \left[- x_2 (x_1^2 + x_2^2) + u_2 - x_1 (x_1^2 + x_2^2) s + u_1 s\right]
+ \left[ (- x_2 - x_1 s) (x_1^2 + x_2^2) s^2\right] \\
&= \frac{d}{dt} X_4 + o(s) \, .
\end{split}
\end{gather*}
\end{example}

\begin{example}
\label{exn4r2PrbNo}
($n = 4$, $m = 2$)
Consider the problem:
\begin{equation*}
\begin{cases}
\dot{x}_1&=u_1(1+x_2)\\
\dot{x}_2&=u_1 x_3\\
\dot{x}_3&=u_2\\
\dot{x}_4&=u_1 x_3^2
\end{cases}
\end{equation*}
with $L=u_1^2+u_2^2$. From Theorem~\ref{r:fundProp}
we get the following necessary conditions for the one-parameter
transformation $h^s = \left(T,X_1,X_2,X_3,X_4,U_1,U_2\right)$
to leave the problem quasi-invariant:
\begin{equation*}
\begin{cases}
\ds \frac{d}{dt} \left.\frac{\partial F}{\partial s}\right|_{s = 0}
     &= \ds 2 u_1 \left.\frac{\partial U_1}{\partial s}\right|_{s = 0}
       + 2 u_2 \left.\frac{\partial U_2}{\partial s}\right|_{s = 0}
       + \left(u_1^2 + u_2^2\right)
       \frac{d}{dt} \left.\frac{\partial T}{\partial s}\right|_{s = 0}\\[0.3cm]
\ds \frac{d}{dt} \left.\frac{\partial X_1}{\partial s}\right|_{s = 0}
     &= \ds  u_1 \left.\frac{\partial X_2}{\partial s}\right|_{s = 0}
     + \left(1 + x_2\right) \left.\frac{\partial U_1}{\partial s}\right|_{s = 0}
     + u_1 \left(1 + x_2\right)
       \frac{d}{dt} \left.\frac{\partial T}{\partial s}\right|_{s = 0}\\[0.3cm]
\ds \frac{d}{dt} \left.\frac{\partial X_2}{\partial s}\right|_{s = 0}
     &= \ds  u_1 \left.\frac{\partial X_3}{\partial s}\right|_{s = 0}
     + x_3 \left.\frac{\partial U_1}{\partial s}\right|_{s = 0}
     + u_1 x_3 \frac{d}{dt} \left.\frac{\partial T}{\partial s}\right|_{s = 0}\\[0.3cm]
\ds \frac{d}{dt} \left.\frac{\partial X_3}{\partial s}\right|_{s = 0}
     &= \ds \left.\frac{\partial U_2}{\partial s}\right|_{s = 0}
     + u_2 \frac{d}{dt} \left.\frac{\partial T}{\partial s}\right|_{s = 0}\\[0.3cm]
\ds \frac{d}{dt} \left.\frac{\partial X_4}{\partial s}\right|_{s = 0}
     &= \ds 2 u_1 x_3 \left.\frac{\partial X_3}{\partial s}\right|_{s = 0}
      + x_3^2 \left.\frac{\partial U_1}{\partial s}\right|_{s = 0}
     + u_1 x_3^2 \frac{d}{dt} \left.\frac{\partial T}{\partial s}\right|_{s = 0} \, .
\end{cases}
\end{equation*}
The conditions are satisfied with $F \equiv 0$ and
\begin{gather*}
\left.\frac{\partial U_1}{\partial s}\right|_{s = 0} = - u_1 \, , \quad
\left.\frac{\partial U_2}{\partial s}\right|_{s = 0} = - u_2 \, , \quad
\frac{d}{dt} \left.\frac{\partial T}{\partial s}\right|_{s = 0} = 2 \, ,\\
\left.\frac{\partial X_1}{\partial s}\right|_{s = 0} = 3 x_1 \, , \quad
\left.\frac{\partial X_2}{\partial s}\right|_{s = 0} = 2 \left(1 + x_2\right) \, , \quad
\left.\frac{\partial X_3}{\partial s}\right|_{s = 0} = x_3 \, , \quad
\left.\frac{\partial X_4}{\partial s}\right|_{s = 0} = 3 x_4 \, .
\end{gather*}
With the transformations
$U_1=u_1(1-s)$, $U_2=u_2(1-s)$, $T=t(1+2s)$,
$X_1=x_1(1+3s)$, $X_2=x_2+2s(1+x_2)$, $X_3=x_3(1+s)$, $X_4=x_4(1+3s)$,
the problem is quasi-invariant:
\begin{gather*}
\begin{split}
L \circ h^s \, \frac{d}{dt} T  &=
\left(u_1^2 + u_2^2\right) + \left(u_1^2 + u_2^2\right) \left(2 s - 3\right) s^2 \, , \\
\varphi_1 \circ h^s \, \frac{d}{dt} T &=
\frac{d}{dt} \left(x_1 (1 + 3 s)\right) - 4 u_1 (1 + x_2) s^3 \, , \\
\varphi_2 \circ h^s \, \frac{d}{dt} T &=
\frac{d}{dt} \left(x_2 + 2s (1 + x_2)\right) - u_1 x_3 (1 + 2s) s^2 \, , \\
\varphi_3 \circ h^s \, \frac{d}{dt} T &=
\frac{d}{dt} \left(x_3 (1 + s)\right) - 2 u_2 s^2 \, , \\
\varphi_4 \circ h^s \, \frac{d}{dt} T &=
\frac{d}{dt} \left(x_4 (1 + 3 s)\right)
+ u_1 x_3^2 \left(1 - 3s - 2 s^2\right) s^2 \, .
\end{split}
\end{gather*}
\end{example}

We will now see how to derive conservation laws from the
knowledge of such $T$, $F$ and $X_i$'s ($i=1,\ldots,n$).


\section{The Noether Theorem and Conservation Laws}
\label{sec:NT}

Now we obtain, as a corollary of Theorem~\ref{r:fundProp},
a far more general Noether theorem
for optimal control problems, which permits
to construct conserved quantities along the Pontryagin
extremals of the problem. Theorem~\ref{r:consLaw}
gives $r$ conservation laws
when problem $(P)$ is quasi-invariant under a family of transformations
containing $r$ parameters.

\begin{theorem}
\label{r:consLaw}
If problem $(P)$ is quasi-invariant under an $r$-parameter family
of transformations \eqref{familyT}
then, for any quadruple
$\left(x(\cdot),u(\cdot),\psi_0,\psi(\cdot)\right)$
satisfying the Pontryagin maximum principle for $(P)$,
the $r$ expressions hold true ($k=1,\ldots,r$):
\begin{multline*}
\psi_0 \left.\frac{\partial F\left(t,x(t),u(t),s\right)}{\partial s_k}
\right|_{s = 0} + \psi(t) \cdot
\left.\frac{\partial X\left(t,x(t),u(t),s\right)}{\partial s_k}\right|_{s = 0} \\
- H(t,x(t),u(t),\psi_0,\psi(t))
\left.\frac{\partial T\left(t,x(t),u(t),s\right)}{\partial s_k}\right|_{s = 0}
\equiv \text{constant} \, ,
\end{multline*}
$t \in [a,b]$, with $H$ the Hamiltonian associated
to the problem $(P)$: $H(t,x,u,\psi_0,\psi)
= \psi_0 L(t,x,u) + \psi \cdot \varphi(t,x,u)$.
\end{theorem}

\begin{remark}
Following the usual terminology (\textrm{cf. e.g.}
\cite[p.~554]{MR1829160}, \cite{MR99b:58095}),
we call a function  $C\left(t,\,x,\,u,\,\psi_0,\,\psi\right)$
which is constant along every Pontryagin extremal
$\left(x(\cdot),\,u(\cdot),\,\psi_0,\,\psi(\cdot)\right)$
of $(P)$,
\begin{equation}
\label{eq:def:CL}
C\left(t,\,x(t),\,u(t),\,\psi_0,\,\psi(t)\right) = k \, ,
\end{equation}
for some constant $k$, a \emph{constant of the motion} or a \emph{first integral}.
The equation \eqref{eq:def:CL} is called the \emph{conservation law} corresponding
to the first integral $C(\cdot,\cdot,\cdot,\cdot,\cdot)$.
\end{remark}

\begin{remark}
As far as everything under consideration, including the Pontryagin
maximum principle, is of a local character, the fact that we restrict
ourserves to state variables in Euclidean spaces $\mathbb{R}^n$
does not lead to any loss of generality. In particular,
Theorem~\ref{r:consLaw} is easily formulated on Manifolds.
\end{remark}

\begin{example}
For the problem considered in Example~\ref{ex:invFO},
we conclude from Theorem~\ref{r:consLaw} that
$2 \psi_0\left(x_1(t) + x_2(t)\right) + \psi_1(t) t + \psi_2(t) t +
\frac{1}{2} \psi_3(t) \left(x_2(t)\right)^2 t$ is constant
along the Pontryagin extremals.
\end{example}

\begin{example}
For the problem in Example~\ref{exn4r2}, the
following first integral follows from Theorem~\ref{r:consLaw}~:
$-\psi_1(t) x_2(t) + \psi_2(t) x_1(t) - \psi_3(t) x_4(t) + \psi_4(t) x_3(t)$.
\end{example}

\begin{example}
\label{ex:CLexn4r2PrbNo}
From Example~\ref{exn4r2PrbNo} and Theorem~\ref{r:consLaw},
the following constant of the motion holds:
\begin{equation}
\label{e:CLexPrbNo}
3 \psi_1(t) x_1(t) + 2 \psi_2(t) (1+x_2(t)) + \psi_3(t) x_3(t)
+ 3 \psi_4(t) x_4(t) -2 t H \, ,
\end{equation}
with
$H = \psi_0 \left((u_1(t))^2 + (u_2(t))^2\right) +
\psi_1(t) u_1(t) \left(1 + x_2(t)\right) + \psi_2(t) u_1(t) x_3(t)
+ \psi_3(t) u_2(t) + \psi_4(t) u_1(t) \left(x_3(t)\right)^2$.
\end{example}

\begin{remark}
All the conservation laws obtained in the previous examples are
not obvious and not expected \textit{a priori}.
However, once obtained, they can easily be
checked, by differentiation, using the corresponding adjoint system
$\dot{\psi} = - \frac{\partial H}{\partial x}$ and the extremality
condition $\frac{\partial H}{\partial u} = 0$. Let us illustrate
this issue for Example~\ref{ex:CLexn4r2PrbNo}. From the adjoint
system we get that $\psi_1$ and $\psi_4$ are constants, while $\psi_2(t)$
and $\psi_3(t)$ satisfy $\dot{\psi_2}(t) =  - \psi_1 u_1(t)$,
$\dot{\psi_3}(t) =  - \psi_2(t) u_1(t) - 2 \psi_4 u_1(t) x_3(t)$.
Having in mind that the problem is autonomous, and therefore the
Hamiltonian $H$ is constant along the extremals
(\textrm{cf.} \cite{torresCM01I14E}),
differentiation of \eqref{e:CLexPrbNo} allow us to write that
\begin{multline*}
3 \psi_1 u_1(t) \left(1 + x_2(t)\right)
- 2 \psi_1 u_1(t) \left(1 + x_2(t)\right) + 2 \psi_2(t) u_1(t) x_3(t)
- \psi_2(t) u_1(t) x_3(t) \\
- 2 \psi_4 u_1(t) \left(x_3(t)\right)^2 +
\psi_3(t) u_2(t) + 3 \psi_4 u_1(t) \left(x_3(t)\right)^2 - 2 H = 0 \, ,
\end{multline*}
that is,
\begin{equation}
\label{eq:CLAftDiff}
\psi_1 \left(1 + x_2(t)\right) u_1(t) + \psi_2(t) x_3(t) u_1(t)
+ \psi_3(t) u_2(t) + \psi_4 \left(x_3(t)\right)^2 u_1(t) = 2 H \, .
\end{equation}
From the definition of the Hamiltonian, equality
\eqref{eq:CLAftDiff} is equivalent to
$H = - \psi_0\left((u_1(t))^2 + (u_2(t))^2\right)$, a relation
that immediately follows from the extremality condition:
\begin{multline*}
\begin{cases}
2 \psi_0 u_1(t) + \psi_1 \left(1 + x_2(t)\right) + \psi_2(t) x_3(t)
+ \psi_4 \left(x_3(t)\right)^2 = 0 \\
2 \psi_0 u_2(t) + \psi_3(t) = 0
\end{cases} \\
\Rightarrow
\begin{cases}
\psi_1 \left(1 + x_2(t)\right) u_1(t) + \psi_2(t) x_3(t) u_1(t)
+ \psi_4 \left(x_3(t)\right)^2 u_1(t) = - 2 \psi_0 \left(u_1(t)\right)^2 \\
\psi_3(t) u_2(t) = - 2 \psi_0 \left(u_2(t)\right)^2 \, .
\end{cases}
\end{multline*}
\end{remark}

\begin{proof} \emph{(Theorem~\ref{r:consLaw})}
Let $\left(x(\cdot),u(\cdot),\psi_0,\psi(\cdot)\right)$ be a Pontryagin
extremal of $(P)$.
Multiplying \eqref{eq:ds0i} by $\psi_0$, \eqref{eq:ds0ii} by $\psi(t)$,
we can write:
\begin{multline}
\label{eq:joined}
\psi_0 \frac{d}{dt} \left.\frac{\partial F}{\partial s_k}\right|_{s = 0}
+ \psi(t) \cdot \frac{d}{dt} \left.\frac{\partial X}{\partial s_k}\right|_{s = 0} \\
= \psi_0 \left(\frac{\partial L}{\partial t}
      \left.\frac{\partial T}{\partial s_k}\right|_{s = 0}
      + \frac{\partial L}{\partial x} \cdot
      \left.\frac{\partial X}{\partial s_k}\right|_{s = 0}
      + \frac{\partial L}{\partial u} \cdot
      \left.\frac{\partial U}{\partial s_k}\right|_{s = 0}
      + L \frac{d}{dt} \left.\frac{\partial T}{\partial s_k}\right|_{s =
      0} \right) \\
+ \psi(t) \cdot \left(\frac{\partial \varphi}{\partial t}
      \left.\frac{\partial T}{\partial s_k}\right|_{s = 0}
      + \frac{\partial \varphi}{\partial x} \cdot
      \left.\frac{\partial X}{\partial s_k}\right|_{s = 0}
      + \frac{\partial \varphi}{\partial u} \cdot
      \left.\frac{\partial U}{\partial s_k}\right|_{s = 0}
      + \varphi \frac{d}{dt} \left.\frac{\partial T}{\partial s_k}\right|_{s =
      0} \right)\, .
\end{multline}
According to the maximality condition of the Pontryagin maximum principle, the function
\begin{equation*}
\psi_0 L\left(t,x(t),U\left(t,x(t),u(t),s\right)\right)
+ \psi(t) \cdot \varphi\left(t,x(t),U\left(t,x(t),u(t),s\right)\right)
\end{equation*}
attains an extremum for $s = 0$.
Therefore for each $k \in \left\{1,\ldots,r\right\}$
\begin{equation*}
\psi_0 \frac{\partial L}{\partial u} \cdot
\left.\frac{\partial U}{\partial s_k}\right|_{s = 0}
+ \psi(t) \cdot \frac{\partial \varphi}{\partial u} \cdot
\left.\frac{\partial U}{\partial s_k}\right|_{s = 0} = 0
\end{equation*}
and \eqref{eq:joined} simplifies to
\begin{multline*}
\psi_0 \left(\frac{\partial L}{\partial t}
      \left.\frac{\partial T}{\partial s_k}\right|_{s = 0}
      + \frac{\partial L}{\partial x} \cdot
      \left.\frac{\partial X}{\partial s_k}\right|_{s = 0}
      + L \frac{d}{dt} \left.\frac{\partial T}{\partial s_k}\right|_{s =
      0} - \frac{d}{dt} \left.\frac{\partial F}{\partial s_k}\right|_{s = 0}\right) \\
+ \psi(t) \cdot \left(\frac{\partial \varphi}{\partial t}
      \left.\frac{\partial T}{\partial s_k}\right|_{s = 0}
      + \frac{\partial \varphi}{\partial x} \cdot
      \left.\frac{\partial X}{\partial s_k}\right|_{s = 0}
      + \varphi \frac{d}{dt} \left.\frac{\partial T}{\partial s_k}\right|_{s =
      0} - \frac{d}{dt} \left.\frac{\partial X}{\partial s_k}\right|_{s = 0}\right)
      = 0 \, .
\end{multline*}
Using the adjoint system
$\dot{\psi} = -\frac{\partial H}{\partial x}$ and the equality
$\frac{dH}{dt} = \frac{\partial H}{\partial t}$ (\textrm{cf.} \cite{torresCM01I14E}),
one easily concludes
that the above equality is equivalent to
\begin{equation*}
\frac{d}{dt} \left(\psi_0 \left.\frac{\partial F}{\partial s_k}\right|_{s = 0}
+ \psi(t) \cdot \left.\frac{\partial X}{\partial s_k}\right|_{s = 0}
- H \left.\frac{\partial T}{\partial s_k}\right|_{s = 0}\right)
= 0 \, .
\end{equation*}
The proof is complete.
\end{proof}


\section{Illustrative Examples}
\label{s:IllExmp}

The following proposition extends the study of the Martinet
flat problem of sub-Riemannian geometry in \cite[\S 4]{delfimEJC}
(see Example~\ref{e:martinetFlat} below) to the general homogeneous case.
\begin{proposition}
\label{r:homogeneousCase}
If there exist constants $\alpha$, $\beta_1, \ldots, \beta_n$,
$\gamma_1, \ldots, \gamma_m$ $\in \mathbb{R}$, such that for all
$\lambda > 0$
\begin{gather}
L\left(\lambda^{\alpha} t,\lambda^{\beta_1} x_1,\ldots,\lambda^{\beta_n} x_n,
\lambda^{\gamma_1} u_1,\ldots,\lambda^{\gamma_m} u_m\right)
= \lambda^{-\alpha} L\left(t,x_1,\ldots,x_n,u_1,\ldots,u_m\right)\, , \label{e:Lhomog} \\
\varphi_i\left(\lambda^{\alpha} t,\lambda^{\beta_1} x_1,\ldots,\lambda^{\beta_n} x_n,
\lambda^{\gamma_1} u_1,\ldots,\lambda^{\gamma_m} u_m\right)
= \lambda^{\beta_i-\alpha}
\varphi_i\left(t,x_1,\ldots,x_n,u_1,\ldots,u_m\right) \, , \label{e:Fihomog} \\
(i=1,\ldots,n) \notag
\end{gather}
then
\begin{equation*}
\sum_{i = 1}^{n} \beta_i \psi_i(t) x_i(t)
- \alpha H\left(t,x(t),u(t),\psi_0,\psi(t)\right) t \equiv \text{constant}
\end{equation*}
along any Pontryagin extremal $\left(x(\cdot),u(\cdot),\psi_0,\psi(t)\right)$
of $(P)$.
\end{proposition}

\begin{proof}
Differentiating \eqref{e:Lhomog} and \eqref{e:Fihomog} with respect to $\lambda$,
and setting $\lambda = 1$, we get
\begin{gather*}
\alpha L(t,x,u) + \alpha \frac{\partial L(t,x,u)}{\partial t} t
+ \sum_{j=1}^{n} \beta_j \frac{\partial L(t,x,u)}{\partial x_j} x_j +
\sum_{k=1}^{m} \gamma_k \frac{\partial L}{\partial u_k} u_k = 0 \, ,\\
\begin{split}
\left(\alpha-\beta_i\right) \varphi_i(t,x,u) +
\alpha \frac{\partial \varphi_i(t,x,u)}{\partial t} t
+ \sum_{j=1}^{n} \beta_j & \frac{\partial \varphi_i(t,x,u)}{\partial x_j} x_j \\
& + \sum_{k=1}^{m} \gamma_k \frac{\partial \varphi_i(t,x,u)}{\partial u_k} u_k = 0 \, .
\end{split}
\end{gather*}
From these equations, one concludes that conditions \eqref{eq:ds0i}
and \eqref{eq:ds0ii} of the fundamental invariance theorem
are fulfilled if we choose
$F \equiv 0$ and a one-parameter family
of transformations satisfying the relations
\begin{equation}
\label{e:cProvProp}
\left.\frac{\partial T}{\partial s}\right|_{s=0} = \alpha t \, ,
\left.\frac{\partial X_i}{\partial s}\right|_{s=0} = \beta_i x_i \, ,
\left.\frac{\partial U_k}{\partial s}\right|_{s=0} = \gamma_k u_k \, .
\end{equation}
For that it suffices to choose $T = \mathrm{e}^{\alpha s} t$,
$X_i = \mathrm{e}^{\beta_i s} x_i \, (i = 1,\ldots,n)$, and
$U_k = \mathrm{e}^{\gamma_k s} u_k\, (k = 1,\ldots,m)$. The problem
is quasi-invariant under these transformations (Definition~\ref{d:invariant})
and the conclusion follows from Theorem~\ref{r:consLaw}.
\end{proof}

\begin{remark}
It is possible to prove the Proposition~\ref{r:homogeneousCase}
with other choices of the parameter family of maps
satisfying  \eqref{e:cProvProp}. For example,
the same conclusion follows from Theorem~\ref{r:consLaw} with
$T = (s+1)^{\alpha} t$,
$X_i = (s + 1)^{\beta_i} x_i$,
$U_k = (s + 1)^{\gamma_k} u_k$,
and $F \equiv 0$, or $T=(1+\alpha s)t$,
$X_i = (1 + \beta_i s) x_i$,
and $U_k = (1 + \gamma_k s) u_k$.
\end{remark}

\begin{example}[$n = 3$, $m = 2$]
\label{e:martinetFlat}
In the Martinet flat problem of sub-Riemannian geometry,
$L = u_1^2 + u_2^2$, $\varphi_1 = u_1$,
$\varphi_2 = u_2$, $\varphi_3 = \frac{u_1 x_2^2}{2}$.
For $\alpha = 2$, $\beta_1 = \beta_2 = 1$, $\beta_3 = 3$,
$\gamma_1 = \gamma_2 = -1$, one concludes from
Proposition~\ref{r:homogeneousCase} that
\begin{equation}
\label{e:CL:martinetFlat}
\psi_1(t) x_1(t) + \psi_2(t) x_2(t) + 3 \psi_3(t) x_3(t)
- 2 H t
\end{equation}
is constant in $t$ along any Pontryagin extremal
\begin{equation*}
\left(x_1(\cdot),x_2(\cdot),x_3(\cdot),u_1(\cdot),u_2(\cdot),\psi_0,
\psi_1(\cdot),\psi_2(\cdot),\psi_3(\cdot)\right)
\end{equation*}
of the problem, with $H$ the Hamiltonian
\begin{equation*}
H(x_2,u_1,u_2,\psi_0,\psi_1,\psi_2,\psi_3)
= \psi_0 \left(u_1^2 + u_2^2\right) + \psi_1 u_1 + \psi_2 u_2
+ \psi_3 \frac{u_1 x_2^2}{2} \, .
\end{equation*}
The first integral \eqref{e:CL:martinetFlat} was first discovered
in \cite{delfim3NCNW}.
\end{example}

Now we will consider optimal control problems subject to
control-affine dynamics with drift. In all the cases our
new version of Noether's theorem is in order.
The application of Theorem~\ref{r:consLaw}
with invariance up to first-order terms of the parameter $s$
will be crucial in the examples,
and therefore the first integrals we obtain can not be deduced
from the previous results in
\cite{MR83k:49054,MR96i:49037,MR1806135,MR1823009,delfimNoether,delfim3NCNW,delfimEJC}.
\begin{example}[$n = 2$, $m = 1$]
Consider problem $(P)$ with $L = u^2$, $\varphi_1 = 1 + y^2$ and
$\varphi_2 = u$~:
\begin{gather*}
\int_a^b \left(u(t)\right)^2 dt \longrightarrow \min \, , \\
\begin{cases}
\dot{x}(t) = 1 + \left(y(t)\right)^2 \, , \\
\dot{y}(t) = u(t)  \, .
\end{cases}
\end{gather*}
From Theorem~\ref{r:fundProp} one gets the following necessary
conditions for the problem to be quasi-invariant under the one-parameter
transformation $h^s = \left(T,X,Y,U\right)$~:
\begin{equation*}
\begin{cases}
\ds \frac{d}{dt} \left.\frac{\partial F}{\partial s}\right|_{s = 0}
     &= \ds 2 u \left.\frac{\partial U}{\partial s}\right|_{s = 0}
          +  u^2 \frac{d}{dt} \left.\frac{\partial T}{\partial s}\right|_{s = 0} \\[0.3cm]
\ds \frac{d}{dt} \left.\frac{\partial X}{\partial s}\right|_{s = 0}
     &= \ds  2 y \left.\frac{\partial Y}{\partial s}\right|_{s = 0} + \left(1 + y^2\right)
     \frac{d}{dt} \left.\frac{\partial T}{\partial s}\right|_{s = 0} \\[0.3cm]
\ds \frac{d}{dt} \left.\frac{\partial Y}{\partial s}\right|_{s = 0}
     &= \ds  \left.\frac{\partial U}{\partial s}\right|_{s = 0}
     + u \frac{d}{dt} \left.\frac{\partial T}{\partial s}\right|_{s = 0} \, .
\end{cases}
\end{equation*}
These conditions are satisfied with $F \equiv 0$, $T = t (1 - 2s)$,
$U = u (1 + s)$, $X = x + 2s (t - 2x)$, and $Y = y (1 - s)$,
for which the problem is quasi-invariant:
\begin{gather*}
L \circ h^s \, \frac{d}{dt} T  = u^2 \left(1 + s\right)^2  \left(1 - 2s\right)
= u^2 - \left(3 + 2s\right) u^2 s^2 = L + o(s) \, , \\
\begin{split}
\varphi_1 \circ h^s \, \frac{d}{dt} T =
\left[1 + y^2 (1 - s)^2 \right] (1 - 2s) & =
\frac{d}{dt} \left[x + 2s (t - 2x)\right] + \left(5y^2 - 2 y^2 s\right) s^2 \\
& = \frac{d}{dt} X  + o(s) \, ,
\end{split} \\
\varphi_2 \circ h^s \, \frac{d}{dt} T = u (1 + s)(1 - 2s)
= u (1-s) - 2 u s^2 = \frac{d}{dt} Y + o(s) \, .
\end{gather*}
From Theorem~\ref{r:consLaw} the following conservation law holds:
\begin{equation*}
2 \psi_x \left(t - 2 x(t)\right) - \psi_{y}(t) y(t) + 2 H t \equiv \text{constant} \, ,
\end{equation*}
where $H = \psi_0 \left(u(t)\right)^2 + \psi_{x}\left[1 + \left(y(t)\right)^2\right]
+ \psi_{y}(t) u(t)$.
\end{example}

In the following two examples we establish conservation laws for the
time-optimal problem.

\begin{example}[$n = 4$, $m = 1$]
\label{ex:TOPn4m1}
Let us consider the minimum-time problem under the control system
\begin{gather*}
\begin{cases}
\dot{x_1}(t) = 1 + x_2(t) \, , \\
\dot{x_2}(t) = x_3(t) \, , \\
\dot{x_3}(t) = u(t) \, , \\
\dot{x_4}(t) = \left(x_3(t)\right)^2 - \left(x_2(t)\right)^2 \, .
\end{cases}
\end{gather*}
In this case the Lagrangian is given by $L = 1$ and in order
to satisfy condition \eqref{eq:ds0i} of the fundamental invariance
theorem we fix $T = t$ (no transformation of the time-variable) and
$F \equiv 0$. The functional is invariant and condition
\eqref{eq:ds0ii} of Theorem~\ref{r:fundProp} simplifies to
\begin{equation*}
\begin{cases}
\ds \frac{d}{dt} \left.\frac{\partial X_1}{\partial s}\right|_{s = 0}
     &= \ds \left.\frac{\partial X_2}{\partial s}\right|_{s = 0} \\[0.3cm]
\ds \frac{d}{dt} \left.\frac{\partial X_2}{\partial s}\right|_{s = 0}
     &= \ds  \left.\frac{\partial X_3}{\partial s}\right|_{s = 0} \\[0.3cm]
\ds \frac{d}{dt} \left.\frac{\partial X_3}{\partial s}\right|_{s = 0}
     &= \ds  \left.\frac{\partial U}{\partial s}\right|_{s = 0} \\[0.3cm]
\ds \frac{d}{dt} \left.\frac{\partial X_4}{\partial s}\right|_{s = 0}
     &= \ds - 2 x_2 \left.\frac{\partial X_2}{\partial s}\right|_{s = 0}
      + 2 x_3 \left.\frac{\partial X_3}{\partial s}\right|_{s = 0} \, .
\end{cases}
\end{equation*}
It is now a simple exercise to conclude that the problem is quasi-invariant,
in the sense of Definition~\ref{d:invariant}, under
$X_1 = (x_1 - t) s + x_1$, $X_2 = x_2 (s + 1)$, $X_3 = x_3 (s + 1)$,
$X_4 = x_4 (2s + 1)$, $U = u (s + 1)$:
\begin{align*}
\frac{d}{dt} X_1 & = \frac{d}{dt} \left[(x_1 - t)s + x_1\right]
= \left(\dot{x_1} - 1\right) s + \dot{x_1} = x_2 s + x_2 + 1 = 1 + X_2 \, , \\
\frac{d}{dt} X_2 & = \frac{d}{dt} \left[x_2 (s + 1)\right]
= \dot{x_2} (s+1) = x_3 (s+1) = X_3 \, , \\
\frac{d}{dt} X_3 & = \frac{d}{dt} \left[x_3 (s + 1)\right]
= u (s+1) = U \, , \\
\frac{d}{dt} X_4 & = \frac{d}{dt} \left[x_4 (2s + 1)\right]
= \dot{x_4} (2s+1) = \left(x_{3}^2 - x_{2}^2\right) (2s+1)
= X_{3}^2 - X_{2}^2 - o(s)\, ,
\end{align*}
with $o(s) = s^2 \left(x_{3}^2 - x_{2}^2\right)$. One obtains
from Theorem~\ref{r:consLaw} the conservation law
\begin{equation}
\label{eq:CLTOPn4m1}
\psi_1(t) \left(x_1(t) - t\right) + \psi_2(t) x_2(t) + \psi_3(t) x_3(t)
 + 2 \psi_4(t) x_4(t) \equiv \text{constant} \, .
\end{equation}
\end{example}

\begin{example}[$n = 3$, $m = 1$]
\label{ex:TOPn3m1}
We consider now the time-optimal problem ($L = 1$) with the control system
\begin{gather*}
\begin{cases}
\dot{x} = 1 + y^2 - z^2 \, , \\
\dot{y} = z \, , \\
\dot{z} = u \, .
\end{cases}
\end{gather*}
From the fundamental invariance theorem one can easily get the
one-parameter transformation
$h^s(t,x,y,z,u) = \left(t,2(x-t)s+x,y(s+1),z(s+1),u(s+1)\right)$,
for which the problem is quasi-invariant ($F \equiv 0$):
\begin{align*}
\frac{d}{dt} X & = \frac{d}{dt} \left[2(x - t)s + x\right]
= (2s+1) (y^2 - z^2) + 1 = 1 + Y^2 - Z^2 - o(s) \, , \\
\frac{d}{dt} Y & = \frac{d}{dt} \left[y (s + 1)\right]
= z (s+1) = Z \, , \\
\frac{d}{dt} Z & = \frac{d}{dt} \left[z (s + 1)\right]
= u (s+1) = U \, ,
\end{align*}
with $o(s) = s^2 \left(y^2 - z^2\right)$. The first integral
associated to the transformation is
\begin{equation}
\label{eq:CLTOPn3m1}
2 \psi_x (x - t) + \psi_y y + \psi_z z \, .
\end{equation}
\end{example}

In Examples~\ref{ex:TOPn4m1} and \ref{ex:TOPn3m1}, if instead of
the time-optimal problem one consider problem $(P)$ with
$J\left[u(\cdot)\right] = \int_a^b u(t) dt \longrightarrow \min$,
the same parameter-transformations are in order choosing appropriate
functions $F$: $F = s x_3$ and $F = s z$ respectively.
The new functionals become invariant up to an exact differential
and the terms $\psi_0 x_3$ and $\psi_0 z$ must be added respectively
to the conservation law \eqref{eq:CLTOPn4m1} and to the constant
of the motion \eqref{eq:CLTOPn3m1}.


\section*{Acknowledgments}

The author acknowledges the financial support of the
program PRODEP III 5.3/C/200.009/2000. The paper has
considerably benefited from the many conversations,
comments and suggestions of Emmanuel Tr\'{e}lat.



\end{document}